\documentclass[12pt]{article}
\usepackage{amscd,amsthm}
\usepackage{latexsym}
\usepackage{amsmath}
\usepackage{amssymb}
\numberwithin{equation}{section}

\def\Re{{\rm Re}}
\def\Im{{\rm Im}}
\def\const{{\rm const}}
\def\Trace{{\rm Trace}}

\def\n{{\cal N}}

\def\p{{\cal P}}

\begin{document}
\title{A Note on Universality of the Distribution of the Largest Eigenvalues 
in  Certain  Sample Covariance Matrices}
\author{Alexander Soshnikov\\University of California\\Department of
Mathematics\\One Shields Avenue\\Davis, CA  95616 USA}
\date{}
\maketitle
\begin{abstract}
Recently  Johansson [19]) and Johnstone ([21]) proved that the distribution
of the (properly rescaled) largest principal component of the complex (real)
Wishart matrix 
$ \ X^*X \ \ (X^tX) \ $
converges to the  Tracy-Widom law as  $n, \ p \ $ (the dimensions of $X$) 
tend to $\infty$ in some ratio $n/p \to \gamma >0.$ 
We extend these results in two directions.  First of all, we prove that the 
joint 
distribution of the first, second, third, etc. eigenvalues of a Wishart 
matrix converges (after a proper rescaling) to the Tracy-Widom distribution.
Second of all, we explain how the  
combinatorial machinery developed  for Wigner random 
matrices in [28]-[30] allows to extend 
the results  by Johansson and Johnstone to the case of $X$ with non-Gaussian 
entries,
provided $ n-p=O(p^{1/3}).$ We also prove that
$ \ \lambda_{max}\leq (n^{1/2} +p^{1/2})^2 +O( p^{1/2} \* \log(p)) 
\  \ (a.e.) \ \ $ for general 
$\gamma >0.$
\end{abstract}

\section{Introduction }

Sample covariance matrices were introduced by statisticians  about seventy 
years ago ([25], [39]). There is a large literature on the subject (see 
e.g. [2]-[6], [9], [12]-[18], [21], [23], [36]-[37]).
We start with the real case.

\subsection{Real Sample Covariance Matrices}

The ensemble consists of $p$-dimensional random  matrices $A_p= X^tX \ $ 
($X^t$ denotes a  transpose matrix), where $X$ is an $n \times p$ matrix with 
independent 
real random entries 
$  x_{ij}, \ 1\leq i \leq n, \ \ 1 \leq j \leq p \ $ such that

\noindent (i)
\begin{align}
{\bf E} x_{ij}& = 0, \\
{\bf E} (x_{ij})^2& = 1, 
\end{align}
$$ 1\leq i \leq n, \ 1\leq j \leq p.$$

To prove the results of Theorems 2, 3  below we will need some 
additional assumptions:

\noindent (ii) The random variables $x_{ij},\ 1\leq i \leq n, 
\ 1\leq j \leq p,$
have symmetric laws of distribution.

\noindent (iii) All moments of these random variables are finite; in 
particular (ii) implies that all odd moments vanish.

\noindent (iv) The distributions of $ x_{ij}, \ $ decay at infinity at least 
as fast as a Gaussian distribution, namely
\begin{equation}
{\bf E} (x_{ij})^{2m}\leq (\const \*  m)^m.
\end{equation}
Here and below we denote by const various positive real numbers that do not 
depend on $n,\ p, \ i,\ j$. 

Complex sample covariance matrices are defined in a similar way.

\subsection{Complex Sample Covariance Matrices}

The ensemble consists of $p$-dimensional random  matrices $A_p= X^*X \ $ 
($X^*$ denotes a complex conjugate matrix), where $X$ is an $n \times p$ 
matrix 
with 
independent complex random entries 
$  x_{ij}, \ 1\leq i \leq n, \ \ 1 \leq j \leq p \ $,  such that

\noindent (i')
\begin{align}
{\bf E} x_{ij}& = 0, \\
{\bf E}(x_{ij})^2 &= 0,\\
{\bf E} |x_{ij}|^2& =1, 
\end{align}
$$ 1\leq i \leq
 n, \ 1\leq j \leq p.$$

The additional assumptions in the complex case mirror those from the real case:

\noindent (ii') The random variables $\Re x_{ij},\  \Im x_{ij},\ 
1\leq i \leq n, \ \ 1 \leq j \leq p$, have symmetric laws of distribution.

\noindent (iii') All moments of these random variables are finite; in particular (ii') implies that all odd moments vanish.

\noindent (iv') The distributions of $\Re x_{ij},\ \Im x_{ij}$ decay at 
infinity at least as fast as a Gaussian distribution, namely
\begin{equation}
{\bf E} |x_{ij}|^{2m}\leq (\const \* m)^m.
\end{equation}
\medskip

\noindent{\bf Remark 1}
The archetypical examples of sample covariance matrices is a $p$ variate
Wishart distribution on $n$ degrees of freedom with identity covariance. 
It corresponds to
\begin{equation}
 x_{ij} \sim N(0,1), \ \  1\leq i \leq n, \ \ 1 \leq j \leq p,
\end{equation}
 in the real 
case, and to
\begin{equation}
\Re x_{ij}, \Im x_{ij} \sim N(0,1), \ \  1\leq i \leq n, 
\ \ 1 \leq j \leq p,
\end{equation}
 in the complex case.
\medskip

It was proved in [23], [17], [37] that if (i) ((i') in the real case) is 
satisfied,
\begin{equation}
n/p \to \gamma \geq 1, \ {\rm as} \ p \to \infty, 
\ \  {\rm and}\ \ 
{ \bf E} |x_{ij}|^{2+\delta} < \const
\end{equation}
then the empirical distribution function of the eigenvalues of $A_p/n$ 
converges
to a non-random limit
\begin{equation}
G_{A_p/n}(x)= \frac{1}{p} \# \{ \lambda^{(p)}_k \leq x, \ k=1, \ldots, n \}
\rightarrow G(x) \ \ (a.s).
\end{equation}
where
$$ \lambda^p_1 \geq \lambda^p_2 \geq \ldots \lambda^p_p $$
are the eigenvalues (all real) of $A_p/n,$
and $G(x)$ is defined by its density $g(x):$
$$g(x)= \begin{cases} \frac{\gamma}{2\pi x}\sqrt{(b-x)(x-a)}, & 
a \leq x \leq b,\\
0, & {\rm otherwise}, \end{cases}$$
$$a= (1-\gamma^{-1/2})^2, \ \ b=(1+\gamma^{-1/2})^2.$$

Since the spectrum of $X \* X^*$ differs from the spectrum of $X^* \* X$ only 
by $ (n-p)$ null eigenvalues, 
the limiting spectral distribution in the case $0<\gamma <1$  remains the same,
except for an atom of mass 
$1-\gamma$ at the origin. From now on we will always assume that $ \ \ p \leq n
,\ \ $ however our results remain valid for $ \ \ p > n \ \ $ as well.

The distribution of the largest eigenvalues attracts a special attention
(see e.g. [21],section 1.2). It was shown by Geman ([15]) in the i.i.d. 
case that if $ {\bf E} |x_{ij}|^{6+\delta} < \infty$ the largest eigenvalue
of $ A_p/n$ converges to $(1+\gamma^{-1/2})^2$ almost surely.
A few years later  Yin, Bai, Krishnaiah and Silverstein ([36], [3])
showed (in the i.i.d. case) that the finiteness of the fourth moment is a 
necessary and sufficient condition for the almost sure convergence (see also  
[27]). These results state that $\lambda_{max}(A_p) = (n^{1/2} + p^{1/2})^2 +
o( n+p).$
However no results were known about the rate of the convergence until recently
Johansson ([19]) and Johnstone ([21]) proved the following theorem in the 
Gaussian (real and complex) cases.
\medskip

 \noindent{\bf Theorem} {\it Suppose that a matrix $A_p=X^t\*X \ $ 
($ A_p=X^* \*X$) has a real (complex) Wishart 
distribution (defined in  Remark 1 above) and $n/p \to \gamma >0.$ Then
$$ \frac{\lambda_{max}(A_p) - \mu_{n,p}}{\sigma_{n,p}}$$
where
\begin{align}
\mu_{n,p}&=(n^{1/2} + p^{1/2})^2,\\
\sigma_{n,p}&=(n^{1/2} + p^{1/2})\*(n^{-1/2}+p^{-1/2})^{1/3}
\end{align}
converges in distribution to the Tracy-Widom law ( $F_1$ in the real case,
$F_2$ in the complex case).}
\medskip

\noindent{\bf Remark 2} Tracy-Widom distribution was discovered by Tracy and 
Widom in [33], [34]. They found that the limiting distribution of the 
(properly rescaled) largest eigenvalue of a Gaussian symmetric
(Gaussian Hermitian) matrix is given by $F_1 (F_2)$, where
\begin{align}
F_1(x)&= \exp \{ -\frac{1}{2}\int_x^{\infty} q(t)+(x-t)\*q^2(t)dt\},\\
F_2(x)&= \exp \{ -\int_x^{\infty} (x-t)\*q^2(t)dt\},
\end{align}
and  $\ q(x) \ $ is such that it solves the Painlev\'e II differential equation
\begin{align}
d^2q(x)/dx^2 &=xq(x)+2q^3(x)\\
q(x) &\sim \textrm{Ai}(x) \ \ {\rm as} \ \  x \to +\infty
\end{align}
where $\ \textrm{Ai}(x) \ $ is the Airy function. Tracy and Widom also derived the 
expressions
for the limiting distribution of the second largest, third largest, etc 
eigenvalues as well.  Since their discovery the field has
exploded with
a number of fascinating papers with applications to combinatorics, 
representation theory, probability, statistics, mathematical physics, 
in which Tracy-Widom  law appears as a 
limiting distribution ( for  recent surveys we refer the reader to 
[1], [10]).
\medskip

\noindent{\bf Remark 3} It should be noted that Johansson studied the complex 
case and Johnstone  did the real case. Johnstone also gave an alternative proof
in the complex case. We also note that Johnstone  has
$n-1$ instead of $n$ in the 
center and scaling constants $\mu_{n,p}, \ \ \sigma_{n,p}$ 
in the real case. While this change clearly does not affect
the limiting distribution of the largest eigenvalues , the choice of $n-1$ is 
more natural if one uses in the proof the asymptotics of Laguerre 
polynomials.
\medskip

\noindent{\bf Remark 4} 
 On a physical level of rigor the results similar to  those from the 
Johansson-Johnstone Theorem (in the complex case) were derived by Forrester 
in [13].
\medskip

While it was not specifically pointed there, the results obtained in [21]
imply that the joint distribution of the first, second, third , 
$\ldots, k$-th,  $\ k=1,2,\ldots \ $ largest eigenvalues converges 
(after the rescaling (1.12),(1.13)) to the limiting distribution derived by 
Tracy-Widom in [33] and [34].
In the complex case one can think about the limiting distribution as the
distribution of the first $k$ (from the right) particles in the determinantal
random point field with the correlation kernel given by the Airy kernel (2.8).
We remind the reader that a random point field is called determinantal
with a correlation kernel $S(x,y)$ if its correlation functions are given by 
\begin{equation}
\rho_k(x_1, \ldots, x_k)= \det_{1\leq i, j\leq k} S(x_i,x_j), \ \ 
k=1,2,\ldots
\end{equation}
(for more information on determinantal random point field we refer the 
reader to [31]). In the real case the situation is slightly more complicated
(correlation functions are given by the square roots of determinants, see 
Section 2, Lemma 1 and Remark 6).
We claim the following result to be true:
\medskip

\noindent{\bf Theorem 1}
{\it The joint distribution of the first, second, third, etc largest 
eigevalues 
(rescaled as in (1.12), (1.13) ) of a real (complex) Wishart matrix 
converges
to the distribution given by the Tracy-Widom law (i.e. the limiting 
distribution of the first, second, etc rescaled eigenvalues for GOE ($ \ \beta
=1, \ $ real case) or  GUE ( $ \beta = 2, \ $ complex case) correspondingly).}
\medskip

Theorem 1 is  proved in Section 2. Our next result generalizes Theorem 1 to 
the non-Gaussian case, provided $ \ n-p = O(n^{1/3}) \ .$
\medskip

\noindent{\bf Theorem 2}
{\it Let a real (complex) sample covariance matrix satisfy the conditions
$ \ ({\rm i}- {\rm iv}) \ (({\rm i'}-{\rm iv'})) \ $  and 
$ \ n-p = O(p^{1/3}) \ .$ Then the joint 
distribution of the first, second, third, etc largest eigenvalues
(rescaled as in (1.12), (1.13)) converge to the Tracy-Widom law with
$ \ \beta = 1 \ (2)\ .$}
\medskip

Similar result for Wigner random matrices was proven in [30].  For other 
results on universality in random matrices we refer the reader to [26],
[11], [7], [20], [8], [22]. 

While we expect the result of Theorem 2 to be true whenever
$ \ n/p \to \gamma >0 \ ,$ we do not know at this moment how to extend our 
technique 
to the case of general $ \ \gamma \ .$  In this paper we  settle for a weaker 
result.

\noindent{\bf Theorem 3}
{\it Let a real (complex) sample covariance matrix satisfy $({\rm i})-
({\rm iv}) \ 
 ( ({\rm i'})-({\rm iv'}))$  and $ \ n/p \to \gamma >0 \ .$
Then 
$$ \ \ a) \ \ {\bf E} \ \Trace \* A_p^{m} =  \frac{(\sqrt{\gamma} +1)
\*\gamma^{1/4}}{2\*
\sqrt{\pi}} \* \frac{ p \*\mu_{n,p}^m}{m^{3/2}}\*(1+o(1))
 \ \ {\rm if \ \ } m=o(\sqrt{p}).$$
$$ \ \ b) \ \ {\bf E} \ \Trace \* A_p^{m} =  
O( \frac{p \* \mu_{n,p}^m}{
m^{3/2}}) \ \ 
 \ \ {\rm if \ \ } m =O(\sqrt{p}). \ \  $$}
\medskip

As a corollary of Theorem 3 we have

\noindent{\bf Corollary 1}
{\it
$$  \lambda_{max}(A_p) \leq \mu_{n,p} +O( p^{1/2} 
\*\log(p) ) \ \ \ 
(a.e.).
$$
}
\medskip

We prove Theorem 1 in Section 2,  Theorem 2 in Section 3 and Theorem 3 and 
Corollary 1 in 
Section 4.

The author would like
to thank Craig Tracy for useful conversations.

\section{ Wishart Distribution}

The analysis in the Gaussian cases is simplified a great deal by the 
exact formulas for the joint distribution of the eigenvalues 
and the $k$-point correlation functions, $k=1,2,\ldots .$
In the complex case the density of the joint distribution of the eigenvalues
is given by ([18]):
\begin{equation}
P_p(x_1, \ldots, x_p)= c_{n,p}^{-1} \prod_{1\leq i < j \leq p}
(x_i -x_j)^2 \prod_{j=1}^{p} x_j^{\alpha_p}\exp(-x_j), \ \ \alpha_p=n-p,
\end{equation}
where $c_{n,p}$ is a normalization constant. Using a standard argument
from Random Matrix Theory ([24] ) one can rewrite 
$P_p(x_1, \ldots, x_p)$ as
\begin{equation}
\frac{1}{p!} \det_{1\leq i,j\leq p} S_p(x_i,x_j)
\end{equation}
where 
\begin{equation}
S_p(x,y)= \sum_{j=0}^{p-1} \varphi^{(\alpha_p)}_j(x)\*
\varphi^{(\alpha_p)}_j(y)
\end{equation}
is the reproducing (Christoffel-Darboux) kernel of the Laguerre 
orthonormalized system
\begin{equation}
\varphi^{(\alpha_p)}_j(x)=\sqrt{\frac{j!}{(j+\alpha_p)!}}\* x^{\alpha_p/2} 
\exp(-x/2)\*L_j^{\alpha_p}(x),
\end{equation}
and $L_j^{\alpha_p}$ are the Laguerre polynomials ([32]).
This allows one to write the $k$-point correlation functions as
\begin{equation}
\rho^{(p)}_k(x_1, \ldots, x_k)= \det_{1\leq i, j\leq k} S_p(x_i,x_j), \ \ 
k=1,2,\ldots,p
\end{equation}
(for more information on  correlation functions we refer the reader to 
[24], [35], [36]).
As a by-product of the  results in [21] Johnstone showed that after the 
rescaling 
\begin{equation}
x= \mu_{n,p} +\sigma_{n,p}\* s
\end{equation}
the (rescaled) kernel 
\begin{equation}
 \sigma_{n,p} S_p(\mu_{n,p} +\sigma_{n,p}\* s_1,
\mu_{n,p} +\sigma_{n,p}\* s_2)
\end{equation}
converges to the Airy kernel
\begin{equation}
S(s_1, s_2)=\frac{A(s_1)\cdot A'(s_2)-A'(s_1)\cdot A(s_2)}{s_1-s_2} =
\ \int_0^{+\infty} \* Ai(s_1 +t) \* Ai(s_2 +t) \* dt.
\end{equation}
The convergence is pointwise and also in the trace norm on any $(t, \infty),
\ \ t\in R^1.$

In the real Wishart case the formula for the joint distribution of the 
eigenvalues was independently discovered by  several groups of 
statisticians at the end of thirties (see [25], [39]):
\begin{equation}
P_p(x_1, \ldots, x_p)= \const_{n,p}^{-1} 
\prod_{1 \leq i < j \leq p}          
|x_i -x_j|  \prod_{j=1}^{p} x_j^{\alpha_p/2}
\exp(-x_j/2), \ \ \alpha_p=n-1-p.
\end{equation}
(note that in the real case $ \ \alpha_p=n-1-p, \ \ $ while in the complex case
it was $ \ n-p. \ )$
The $k$-point correlation function has a form similar to (2.2), (2.3) however
it is  now equal to a square root of the determinant, and 
$K_p(x,y)$ is a  $2 \times 2$ matrix kernel 
(see e.g. [38],[21]):
\begin{equation}
\rho^{(p)}_k(x_1, \ldots, x_k)= \bigl(\det_{1\leq i, j\leq k} K_p(x_i,x_j)
\ \ 
\bigr)^{1/2},\ \ \ k=1,\ldots, p,
\end{equation}
where (in the even $p$ case)
\begin{align}
K^{(1,1)}_p (x,y)&= S_p(x,y) + \psi(x)\* (\epsilon\*\phi)(y)\\
K^{(1,2)}_p (x,y)&= (S_p\*D)(x,y) -\psi(x)\* \phi(y)\\
K^{(2,1)}_p (x,y)&= (\epsilon\*S_p)(x,y) - \epsilon (x-y) + 
(\epsilon \* \psi )(x)\* (\epsilon\phi)(y)\\
K^{(2,2)}_p (x,y)&= K^{(1,1)}_p (y,x), 
\end{align}
operator $\epsilon \ $ denotes convolution with the kernel 
$$ \ \epsilon(x-y) \ =
\frac{1}{2} {\rm sign} (x-y), \ \ \ (S\*D)(x,y)= - \frac{\partial S(x,y)}{
\partial y}, $$ and $\ \psi (x), \ \ \phi (x) \ \ $ are defined as follows
\begin{align}
& \psi(x)= (-1)^p  \frac{(p\*(p+\alpha_p))^{1/4}}{2^{1/2}} \bigl( \sqrt{p+\alpha_p}
\* \xi_p (x) - \sqrt{p} \* \xi_{p-1}(x) \\
& \phi(x)= (-1)^p \* \frac{(p\*(p+\alpha_p))^{1/4}}{2^{1/2}} \bigl( \sqrt{p}
\* \xi_p (x) - \sqrt{p+\alpha_p} \* \xi_{p-1}(x) \\
& \xi_k (x)= \varphi_k^{(\alpha_p)}(x)/x.
\end{align}

\medskip
\noindent{ \bf Remark 5}
The formulas for $ \ K_p(x,y) \ $ in the odd $p$ case are slightly different.
However since we are interested in the asymptotic behavior of the largest 
eigenvalues it is enough to consider only even $p$ case.
Indeed, one can carry very similar calculations in the 
odd $p$
case and obtain the same limiting kernel $ \ \ K(x,y) \ \ $ as we got in 
Lemma 1. 
Or one may observe
that the limiting distribution of the largest (rescaled)
eigenvalues must be the same in the even $p$ and odd $p$ cases as implied
by the following argument. Consider  an $ \ \ (n+p)\times(n+p) \ \ $
real symmetric (self-adjoint) matrix $B=(b_{ij}), \ \ 1 \leq i,j \leq n+p,$ 
$$ b_{ij}= \begin{cases}
           x_{i,j-n}, &{\rm if}  \ \ 1 \leq i \leq n, \ \ n+1 \leq j \leq n+p\\
          \bar x_{j, i-p}, & {\rm if}\ \  p+1 \leq i \leq n+p, \ \ 1\leq j \leq
                   p\\
                   0, &{\rm otherwise}, \end{cases}$$
Then the non-zero eigenvalues of $\ B^2 \ $ and $ \ X^*\*X \ $ coincide.
If we now consider a matrix $ \ \tilde X \ $ obtained by deleting the first row
and the last column of $ \ X \ $ and construct the corresponding matrix
$ \ \tilde B \ $, then by the mini-max principle we have 
$ \ \ \lambda_k (B) \geq
\lambda_k(\tilde B), \ \ k=1,2,\ldots .$ Repeating this procedure once more
we see that the $k$-th eigenvalue of $\ X^*\*X \ $ for odd $p$ is sandwiched
between the $k$-th eigenvalues for $\ p+1 \ $ and $ \ p-1 \ .$
\medskip

The machinery developed in [21] allows us to obtain the following result 
about the pointwise convergence of the entries of $ \ K_p (x,y). \ $
\medskip

\noindent{ \bf Lemma 1}
{\it 
\begin{align}
a) \ \ & \sigma_{n,p} \* K^{(1,1)}_p(\mu_{n,p} +\sigma_{n,p}\* s_1,
\mu_{n,p} +\sigma_{n,p}\* s_2) \to S(s_1, s_2) + \frac{1}{2} \* \textrm{Ai}
(s_1) \* 
\int_{-\infty}^{s_2} \textrm{Ai}(t)\* dt,\\
& \sigma_{n,p} \* K^{(2,2)}_p(\mu_{n,p} +\sigma_{n,p}\* s_1,
\mu_{n,p} +\sigma_{n,p}\* s_2) \to S(s_2, s_1) + \frac{1}{2} \* \textrm{Ai}
(s_2) \* 
\int_{-\infty}^{s_1} \textrm{Ai}(t)\* dt,\\
& \sigma_{n,p}^2 \* K^{(1,2)}_p(\mu_{n,p} +\sigma_{n,p}\* s_1,
\mu_{n,p} +\sigma_{n,p}\* s_2) \to -\frac{1}{2} \* \textrm{Ai}(s_1) \* 
\textrm{Ai}(s_2) -
\frac{\partial}{\partial \* s_2} \* S(s_1, s_2) ,
\end{align}
$$
 K^{(2,1)}_p(\mu_{n,p} +\sigma_{n,p}\* s_1,
\mu_{n,p} +\sigma_{n,p}\* s_2) \* \to  \* - \int_0^{+\infty} \* du \* 
\bigl(\int_{s_1 +u}^{+\infty} \textrm{Ai}(v) dv \bigr) \* \textrm{Ai}(s_2 +u) 
 $$
\begin{equation}
-\epsilon(x-y) \ +\frac{1}{2}\* \int_{s_2}^{s_1} \* Ai(u)\*du \ 
+\frac{1}{2} \* \int_{s_1}^{+\infty} \* \textrm{Ai}(u) \* du
\* \int _{-\infty}^{s_2} \* \textrm{Ai}(v) \* dv.
\end{equation}
b) Convergence in (2.18)-(2.21) is uniform on $\ \ [\tilde s_1, +\infty) \times
[\tilde s_2, +\infty) \ \ $ as $ \ \ p \to \infty \ \ $for any 
$ \ \ \tilde s_1 > -\infty, 
\ \ \tilde s_2 >-\infty \ \ . $
It is also true that the error terms are $ \ \ O(e^{-const \*(s_1 +s_2)})\ \ $ 
uniformly in
$ \ p \ $ with some
$ \ \ const >0.$ 
}
\medskip

\noindent{\bf Remark 6}

Lemma 1 implies the convergence of the rescaled $k$-point correlation functions
$ \ \ \sigma_{n,p}^k \*\rho^{(p)}_k(x_1, \ldots, x_k), \ \   
\ x_i=\mu_{n,p}+\sigma_{n,p}\* s_i,\ \ i=1,\ldots,k,
\ \   k=1,2,\ldots \ \ $ to
$$
\rho_k(s_1, \ldots, s_k)= \bigl(\det_{1\leq i, j\leq k} K(s_i,s_j), 
\ \ 
\bigr)^{1/2},
$$
where the entries of $ \ \ K(s,t)=\bigl ( K_{ij}(s,t) \bigr )_{i,j=1,2} \ \ $
are given by the r.h.s. of (2.18)-(2.21). The limiting correlation functions 
coincide with the limiting correlation functions at the edge of the spectrum
in the Gaussian Orthogonal Ensemble (see e.g. [14]) (it also should be noted
that the formulas (1.15)-(1.16) we gave in [30] for $ \ \ K(s,t) \ \ $
must be replaced by (2.18)-(2.21)).
\medskip   

\noindent{\bf Proof of Lemma 1}  

The proof is a consequence of (2.11)-(2.14), (1.12)-(1.13) and the asymptotic 
formulas for the Laguerre polynomials $ \ \ L_j^{\alpha_p}(x) \ \ ,
\alpha_p \to \infty, \ \ j \sim \alpha_p, \ \ $ near the 
turning point derived in [21].
Below we prove (2.18) and (2.21).
(2.19) immediately follows from (2.14) and (2.18). 
(2.20) is established in a similar way to (2.18), (2.21).
To prove (2.18) we employ a very useful integral representation
for $ \ \ S_p(x,y) \ \ $  ([38]):
\begin{equation}
S_p(x,y)= \int_0^{+\infty} \ \phi(x+z)\*\psi(y+z) + \psi(x+z)\*\phi(y+z) \* dz,
\end{equation}
where $ \ \ \phi(x), \ \psi(x) \ \ $ are defined in (2.15)-(2.17).

The asymptotic behavior of $ \ \ \phi(x), \ \ \psi(x) \ \ $ was studied by
Johnstone ([21]) who proved 
\begin{equation}
\sigma_{n,p}\* \phi(\mu_{n,p} + \sigma_{n,p}\* s), \ \ 
\sigma_{n,p}\* \psi(\mu_{n,p} + \sigma_{n,p}\* s) \ \to \frac{1}{\sqrt{2}}\*
\textrm{Ai}(s)
\end{equation}
and that the l.h.s. at (2.24) is exponentially small for large 
$\ s_1, \ s_2 \ $
(uniformly in $p$.) While Johnstone stated only pointwise convergence in 
(2.22) 
his results (see (3.7), (5.1), (5.19), (5.18), (5.22)-(5.24) and (6.11) 
from [21])
actually imply that the convergence is uniform on any $ \ [s, +\infty) .$
This together with (2.22) gives us
\begin{equation}
\sigma_{n,p}\*S_p(\mu_{n,p} + \sigma_{n,p}\* s_1, 
\mu_{n,p} + \sigma_{n,p}\* s_2) \to S(s_1, s_2),
\end{equation}
where the convergence is uniform on any $ \ [\tilde s_1, \infty) \times
[\tilde s_2, \infty) \ .$
To deal with the second term at the r.h.s. of (2.11),
\begin{equation}
 \ \psi(x)\* (\epsilon\*\phi)(y) =\psi(x) \* \bigl ( \frac{1}{2} \*
\int_0^{\infty}\ \phi(u) \* du -\int_y^{\infty}\ \phi(u) \* du \bigr ),
\end{equation}
we use
\begin{equation}
\frac{1}{2} \*\int_0^{\infty}\ \phi(u) \* du \to\frac{1}{\sqrt{2}}
\end{equation}
(see [21], Appendix A7).
(2.23), (2.25)-(2.26) imply
\begin{equation}
\sigma_{n,p} \* \psi(\mu_{n,p}+\sigma_{n,p}\* s_1)\* (\epsilon\*\phi)
(\mu_{n,p}+\sigma_{n,p}\* s_2) \to 
\frac{1}{2} \* \textrm{Ai}(s_1) \* \int_{-\infty}^{s_2} \textrm{Ai}(t)\* dt.
\end{equation}
This proves (2.18).

To establish (2.21) we consider separately
$ \ \  (\epsilon\*S_p)(x,y) \ \ $ and \\$ 
 (\epsilon \* \psi )(x) \* (\epsilon\phi)(y). \ \ $
We have

\begin{align}
\epsilon\*S_p(x,y) &= \bigl ( \frac{1}{2} \* \int_0^{+\infty} \*  du -
\int_x^{+\infty} \* du \bigr ) \* 
\int_0^{+\infty} \* \phi(u+z)\*\psi(y+z) + \psi(u+z)\*\phi(y+z) \* dz
\\
& =\frac{1}{2}\* \int_0^{+\infty} \* \bigl ( \int_z^{+\infty} \* \phi(u) \* du
\* \psi(y+z) \bigr ) dz \\
& - \int_0^{+\infty} \* \bigl ( \int_{x+z}^{+\infty} \* \phi(u) \* du
\* \psi(y+z) \bigr ) dz \\
& +\frac{1}{2}\* \int_0^{+\infty} \* \bigl ( \int_z^{+\infty} \* \psi(u) \* du
\* \phi(y+z) \bigr ) dz \\
& - \int_0^{+\infty} \* \bigl ( \int_{x+z}^{+\infty} \* \psi(u) \* du
\* \phi(y+z) \bigr ) dz 
\end{align}

Let us fix $ \ \ s_1, \ s_2$ and consider 
\begin{equation}
x= \mu_{n,p} + \sigma_{n,p}\* s_1, \ \ y= \mu_{n,p} +\sigma_{n,p} \*s_2.
\end{equation}
It follows from (2.23) that  
the integrals (2.30) and (2.32)  converge to
$$ - \frac{1}{2} \* \int_0^{+\infty} \* du \* 
\bigl(\int_{s_1 +u}^{+\infty} \textrm{Ai}(v) dv \bigr) \* \textrm{Ai}
(s_2 +u) .$$
Let us now write (2.29) as
\begin{align}
& \frac{1}{2}\* \int_0^{+\infty} \* \bigl ( \int_0^{+\infty} \* \phi(u) \* du
\* \psi(y+z) \bigr ) dz \\
& -\frac{1}{2}\* \int_0^{+\infty} \* \bigl ( \int_0^{z} \* \phi(u) \* du
\* \psi(y+z) \bigr ) dz \\
& =\frac{1}{\sqrt{2}}\* \int_0^{+\infty} \* 
\psi(y+z)  dz \\
& -\frac{1}{2}\* \int_0^{+\infty} \* \bigl ( \int_0^{z} \* \phi(u) \* du
\* \psi(y+z) \bigr ) dz 
\end{align}
Using (2.23) one can see that  (2.36) 
converges to
$ \frac{1}{2}\* \int_{s_2}^{\infty} \* \textrm{Ai}(u)\*du. $\\
The integral (2.37) tends to zero as $ \ p \to \infty \ $. Indeed, 
suppose  that $ \ n-p \to +\infty \  \ $(the case $ \ \ n-p =O(1) \ \ $ can be 
treated  by using  the classical asymptotic formulas  for 
Laguerre polynomials  for 
fixed $ \ \alpha \ $ (see e.g.
[32])).  Let us  write
$\ \  \int_0^{+\infty} \* \bigl ( \int_0^{z} \* \phi(u) \* du
\* \psi(y+z) \bigr ) dz \ \ $as
\begin{equation}
\int_0^{\sqrt{p}} \* \bigl ( \int_0^{z} \* \phi(u) \* du
\* \psi(y+z) \bigr ) dz
 +\int_{\sqrt{p}}^{\infty} \* \bigl ( \int_0^{z} \* \phi(u) \* du
\* \psi(y+z) \bigr ) dz.
\end{equation}

Similar calculations to the ones from Appendix 7 of [21]  show that
for $ \ \ z<\sqrt{p}\ \ $
$$ \int_0^{z} \* \phi(u) \* du= O( (const \  p)^{-(n-p)/4}), \ \ 
{\rm where} \ \ const >0.$$
This estimate coupled with the following   (rather rough) bounds
$$
 \int_0^{\sqrt{p}} \* 
|\psi(y+z)| dz \leq p^{1/4} (\int_y^{\infty} \* \psi(z)^2 \* dz)^{1/2} \leq$$
$$ const \  p^{1/4} \* \bigl (
(\int_y^{\infty} \* \varphi_p^{\alpha_p}(z)^2 \* dz)^{1/2}
+ (\int_y^{\infty} \* \varphi_{p-1}^{\alpha_p}(z)^2 \* dz)^{1/2} \bigr )\\
=O(p^{1/4})
$$
take care of the first term in (2.38).
If $ \ \ z \geq \sqrt{p} \ \ $ one has
$$ |\psi(y+z)| = |\psi( \mu_{n,p} +\sigma_{n,p} \* ( s_2 + z/\sigma_{n,p}))| <
\exp( -const \* (s_2 + z/p^{1/3})), \ \ \ \ const >0,$$
where we have used   the exponential  decay of 
$ \ \ \psi(\mu_{n,p} +\sigma_{n,p} \* s) \ \ $ for large $ \ s \ $ 
(see (2.23), (2.15)-(2.17) 
 and [21], formula (5.1) ). Since
$$| \int_0^z \* \phi(u) \* du| \leq \sqrt{z} \* 
\bigl (\int_0^z \* \phi(u)^2 \* du \bigr )^{1/2}\\
\leq const \  p \*\sqrt{z},$$
we conclude that (2.37) is o(1). Using 
$ \ \ \ \int_0^{\infty} \* \psi(u) \* du=0 \ \ $ one can prove in a similar 
fashion that (2.31) is also o(1). To establish (2.21) we are left with 
estimating   $$ \ \ (\epsilon \* \psi )(x) \* (\epsilon\phi)(y) =
\bigl ( \frac{1}{2} \* \int_0^{\infty} \psi(u) \* du - \int_x^{\infty}
\* \psi(u) \* du \bigr ) \* \bigl ( \frac{1}{2} \* \int_0^{\infty} \phi(u) 
\* du - \int_y^{\infty} \* \phi(v) \* dv \bigr )$$
$$= \bigl ( -\int_x^{\infty} \* \psi(u) \* du \bigr ) 
\bigl ( \frac{1}{\sqrt{2}} +o(1) -
\* \int_y^{\infty} \* \phi(v) \* dv \bigr ). $$
Using (2.23) and (2.33) we derive that the last expression converges to\\
$ \ \ 
-\frac{1}{2}\* \int_{s_1}^{\infty} \* \textrm{Ai}(u)\*du \ 
+\frac{1}{2} \* \int_{s_1}^{+\infty} \* \textrm{Ai}(u) \* du
\* \int _{-\infty}^{s_2} \* \textrm{Ai}(v) \* dv. \ \ $ This finishes the 
proof of 
(2.21). To obtain (2.20) we use (2.23) and
\begin{equation}
\sigma_{n,p}^2\* \phi'(\mu_{n,p} + \sigma_{n,p}\* s), \ \ 
\sigma_{n,p}^2\* \psi'(\mu_{n,p} + \sigma_{n,p}\* s) \ 
\to \frac{1}{\sqrt{2}}\*
\textrm{Ai}'(s).
\end{equation}
which follows from the machinery developed in [21].
Lemma 1 is proven.
\medskip

Theorem 1
now follows from
\medskip

\noindent{\bf Lemma 2}
{\it Suppose that we are given random point fields $ \ {\bf F}, \  
\ {\bf F_n}, \ \ n=1,2,\ldots$ with the $k$-point correlation functions
$\rho_k(x_1, \ldots, x_k), \ 
\ \rho^{(n)}_k(x_1, \ldots, x_k) \ \ \ k=1,2,\ldots \ $ such that
the number of particles in $ \ (a, \infty) \ $ (denoted by
$\ \#(a,\infty) \ $) is finite ${\bf F}-$a.e. for any $ \ \ a>-\infty \ \ $
and $ \ {\bf F } \ $ is uniquely determined by its correlation functions.
Then  the following diagram holds:  $$ \ \ \ d) \Longrightarrow \ c) 
\Longrightarrow \ b) \iff \ a),$$
where

a) The joint distribution of the first, second, $\ldots, k-$th rightmost 
particles
in ${\bf F_n}$ converges to the joint distribution of the first, second, 
$\ldots, k-$th rightmost particles in ${\bf F}$ for any $k \geq 1.$

b) The joint distribution of $ \ \#(a_1, b_1), \ldots, \#(a_l, b_l), \ \ 
l\geq 1 \ $ in ${\bf F_n}$  converges to the corresponding distribution in
${\bf F}$ for any collection of disjoint intervals
$\ \ (a_1, b_1), \ldots, (a_l, b_l), \ \ 
a_j > -\infty, b_j \leq +\infty, j=1,\ldots, l, \ \  l=1, \ldots.$ 

c) $ \ \ \ \rho_k(x_1, \ldots, x_k)$ is integrable on $[t, \infty)^k$ for any
$t \in R^1, \ \  k=1,2,\ldots$
and
\begin{align}
&\int_{(a_1, b_1)^{k_1}\times \ldots \times (a_l, b_l)^{k_l}}
 \rho_k^{(n)}(x_1,\ldots, x_k) dx_1 \ldots dx_k \to \\
&\int_{(a_1, b_1)^{k_1}\times \ldots \times (a_l, b_l)^{k_l}} \rho_k(x_1,\ldots, x_k) dx_1 \ldots dx_k 
\end{align}
for any disjoint intervals $\ \ (a_1, b_1), \ldots, (a_l, b_l), \ \ 
a_j > -\infty, b_j \leq +\infty, j=1,\ldots, l, \ \  l=1, \ldots, k, 
\ \ k_1 + \ldots + k_l= k, 
\ \ k=1,2,\ldots.$

d) For any $ \ k \geq 1 \ $ the Laplace transform
$$ \int \* \exp( \sum_{j=1,\ldots, k} t_j\*x_j) \*
\rho_k(x_1,\ldots, x_k) dx_1 \ldots dx_k  $$
is finite for  $ t_1 \in [c_1^{(k)},d_1^{(k)}] \ \ldots, t_k \in
[c_k^{(k)}, d_k^{(k)}], \ \ $ where $ \ \ c_j^{(k)}< d_j^{(k)}, 
\ \ d_j^{(k)} >0, \ \ 
j=1, \ldots, k, $ and
\begin{align}
&\int  \* \exp( \sum_{j=1,\ldots, k} t_j\*x_j) \*
\rho_k^{(n)}(x_1,\ldots, x_k) dx_1 \ldots dx_k  \ \ \to \\
&\int \* \exp( \sum_{j=1,\ldots, k} t_j\*x_j) \*
\rho_k(x_1,\ldots, x_k) dx_1 \ldots dx_k  
\end{align}
for such $ t_1, \ \ldots, t_k  \ $ as $ \ n \to \infty$.}
\medskip

\noindent{\bf Proof of Lemma 2} $ \ \ d) \Longrightarrow \ c) $

Suppose that d) holds. Fix some positive $ \ \tilde t_1 \in (c_1^{(k)}, 
d_1^{(k)}),\ldots 
\tilde t_k \in (c_k^{(k)}, d_k^{(k)}) \ .$ Denote by
$H_n (dx_1, \ldots, dx_k), \ \ H(dx_1, \ldots, dx_k),$  the probability 
measures on $R^k $ with the densities
$$ h_n(x_1, \ldots, x_k)= Z_n^{-1} \* \exp( \sum_{j=1,\ldots, k} 
\tilde t_j\*x_j) \* \rho_k^{(n)}(x_1,\ldots, x_k), $$
$$ h(x_1, \ldots, x_k)= Z^{-1} \* \exp( \sum_{j=1,\ldots, k} 
\tilde t_j\*x_j) \* \rho_k(x_1,\ldots, x_k), $$
where $ Z_n, \ \ Z \ \ $ are the normalization constants (it is easy to see 
that $ \ \ Z_n \to Z. \ \ $).
The constructed sequence of probability measures is tight (by Helly 
theorem), moreover
their distributions   decay (at least) exponentially for large 
(positive and negative) $ \ x_1, \ldots,
x_k \ \ $
uniformly in $n$. It follows from the tightness of $ \{H_n \}$ that
all  we have to show is that any limiting point of $H_n$ coincides with
$H.$ Suppose that a subsequence of $H_n$ weakly converges to $\bar H.\ $
Then $\bar H$ must have a finite Laplace transform
for $ \ \ c_1^{(k)}-\tilde t_1 \leq Re \ t_1 \leq d_1^{(k)}-\tilde t_1, 
\ldots, c_k^{(k)} -\tilde t_k \leq Re \  t_k \leq d_k^{(k)}-\tilde t_k \ \ $ 
and 
the Laplace 
transforms of $H_n$ must converge to the Laplace transforms of $\ \bar H \ $
in this strip. Since the Laplace transforms of $ \ \bar H, \ H \ $ are 
analytic in the strip and coincide for
$ \ \ t_1 \in [c_1^{(k)},d_1^{(k)}] \ \ldots, t_k \in [c_k^{(k)},d_k^{(k)}] 
\ \ $ they must coincide in the whole strip. 
Applying the inverse Laplace transform we obtain that $\bar H$ coincides with
$H.$ It follows then that
$$\int_{(a_1, b_1)\times \ldots \times (a_k, b_k)}
 \rho_k^{(n)}(x_1,\ldots, x_k) dx_1 \ldots dx_k \to $$
$$
\int_{(a_1, b_1)\times \ldots \times (a_k, b_k)} \rho_k(x_1,\ldots, x_k) 
dx_1 \ldots dx_k 
$$
for any finite $ \ \ a_j < b_j, \ \ j=1, \ldots, k, \ \ $ and the exponential 
decay of  $ \ \ \rho_k^{(n)}(x_1,\ldots, x_k), \ \ \rho_k(x_1,\ldots, x_k), 
\ \ $ for large positive $ \ \ x_1, \ldots, x_k, \ \ $ implies that this
still holds for $ \ \ b_j=+\infty, \ \ j=1,\ldots, k.$

$ \ \ c) \Longrightarrow b) $

We remind the reader that the integral in (2.40) is equal
to the  $ \ (k_1, \ldots, k_l)$-th factorial moment
$$ {\bf E} \prod_{j=1,\ldots,l} \frac{ \bigl (\# (a_j, b_j) \bigr )!}{ \bigl (
\# (a_j, b_j) - k_j \bigr )!} $$
of the numbers of particles in the disjoint intervals
$ \ (a_1, b_1), \ldots, (a_l, b_l). \ \ $
Since the probability distribution of the random point field {\bf F} is 
uniquely determined by its correlation functions, the joint distribution of 
the numbers of particles in the boxes is uniquely determined by the moments, 
and therefore the convergence of moments implies the convergence of the 
distributions of $ \ \#(a_1, b_1), \ldots, \#(a_l, b_l). \ \ $ 

$ \ \ b) \iff a) $

Trivial. Observe that
$$ P \bigl ( \lambda_1 \leq  s_1, \lambda_2 \leq s_2, \ldots, 
\lambda_k \leq s_k \bigr )= $$
$$ P \bigl ( \# (s_1, +\infty)=0, \# (s_2, +\infty) \leq 1, \ldots, 
\# (s_k, +\infty) \leq k-1 \bigr ).$$
Lemma 2 is proven.
\medskip

\noindent{\bf Proof of Theorem 1}
It is worth noting that the limiting random point fields 
defined in (1.18), (2.8) (complex case) and in Remark 6 
(real case) are uniquely determined by their correlation functions (see e.g.
[31]).
To prove Theorem 1 in the complex case we use  a general fact for the 
ensembles with determinantal correlation functions
that the generating function of the numbers of particles in the boxes 
is given by the Fredholm determinant 
\begin{equation}
{\bf E} \ \prod_{j=1,\ldots, l} \* z_j^{\#(a_j,b_j)}  = \det (Id + 
\sum_{j=1,\ldots, k} (z_j-1) \* S_p \* \chi_{(a_j, b_j)}
\end{equation}
(see e.g. [35], [31]), where $ \ \chi_{(a, b)} \ $ is the operator of the
multiplication by the indicator of $\ (a,b) \ .$ Trace class convergence
of $S_p$ to $K$ on any $ \ (a, \infty), \ a> -\infty, \ $ implies the 
convergence of the Fredholm determinants, which together
with Lemma 2 proves Theorem 1 in the complex case. To prove Theorem 1 in the 
real case we observe that Lemma 1 implies that after rescaling 
$ \ \ x_i= \mu_{n,p}+\sigma_{n,p}\* s_i, \ \ \ i=1,2,\ldots \ $ condition 
(2.40)-(2.41) of Lemma 2,
part c) is satisfied.
Theorem 1 is proven.
\medskip

\section{ Proof of Theorem 2}

The proof of Theorem 2 heavily relies on the results obtained in [28]-
[30]. We start with
\medskip

\noindent{\bf Lemma 3}
{\it 
Let $A_p$ be either a real sample covariance matrix  $({\rm i})-({\rm iv})$
or complex sample covariance matrix 
$ \ \bigl (({\rm i'})-({\rm iv'})\bigr ) \ $ and $ \ n-p = O(p^{1/3}) \ $
as $\ p \to \infty.$
Then there exists some  $ \ const >0 \ $ such that for any 
$t_1,\ t_2,\dots t_k>0$ and
$$m^{(1)}_p =[t_1\cdot p^{\frac{2}{3}}],
\ldots ,m^{(k)}_p = [t_k\cdot p^{\frac{2}{3}}],$$
the following estimate holds:
\begin{enumerate}
\item[a)]
\begin{equation}
{\bf E}\prod^k_{i=1}\ \Trace\ A_p^{m^{(i)}_p}\leq const^k \* 
\prod^k_{i=1} \ \frac{\mu_{n,p}^{ m_p^{(i)}}}{
t_i^{3k/2}}\*\exp (const \ \sum^k_{i=1}t^3_i)
\end{equation}

\item[b)] If  $ \ A_p, \ \tilde A_p \ $ 
belong to two different ensembles of random real
(complex) sample covariance matrices satisfying $({\rm i})-({\rm iv})\  \bigl 
( ({\rm i'})-({\rm iv'})
\bigr ),$ and $ \ \ n-p=O(p^{1/3}),$
then 
\begin{equation}
{\bf E} \prod^k_{i=1}\ \Trace\ A_p^{m^{(i)}_p} -
{\bf E} \prod^k_{i=1}\ \Trace\ \tilde A_p^{m^{(i)}_p} 
\end{equation}
tends to zero as $\ \ p \to \infty.$

\end{enumerate}
}
\medskip

\noindent{\bf Proof of Lemma 3}

Lemma 3 is the analogue of Theorem 3 in [30] and is proved in the same way.
Since the real and the complex cases are very similar, we will consider here
only the real case. As we explained earlier, we can assume without loss 
of 
generality that $ \ \ p \leq n \ \ .$
Our arguments will be the most transparent when $ \ k=1 \ $  and 
the matrix entries
$ \ \ \{ x_{ij} \}, \ \ 1\leq i \leq n, \ \ 1 \leq j \leq p \ $ are 
identically distributed.  Construct a  $ \ n\times n \ \ $ 
random real symmetric Wigner matrix $ \ M_n=(y_{ij}), \ 1 \leq i,j \leq n \  $
such that $ \ y_{ij}= y_{ji}, \ \ i \leq j $ are independent identically 
distributed random variables with the same distribution as $ \ x_{11} \ .$
Then
\begin{equation}
{\bf E} \Trace\ A_p^{m_p} \leq {\bf E} \Trace \ M_n^{2\*m_p}.
\end{equation}
To see this we consider separately the left and the right hand sides of the 
inequality. We start by calculating the mathematical expectation of \\ $ 
\Trace  \ A^{m_p}_p$.
Clearly,
\begin{equation}
{\bf E}
 \ \Trace\ A^{m_p}_p=  \sum_{\p} {\bf E}\ x_{i_1,i_0}
\* x_{i_1, i_2}\* x_{i_3, i_2}\* x_{i_3,i_4}
\dots x_{i_{2\* m_p -1}, i_{2 \* m_p -2}}
\* x_{i_{2 \* m_p-1, i_0}}.
\end{equation}
The sum in (3.4) is taken over all closed paths 
$\p =\{i_0,i_1,\ldots ,i_{2\* m_p -1},i_0 \}$, with a distinguished origin, 
in the set $\{1,2,\dots n\}$ with the condition
$$ \ {\bf C1. \ \ }\ \ i_t \in \{1,2,\dots p\} \ \ {\rm for \ \ odd \ \ } t$$
satisfied.
  We consider the set of vertices 
$\{1,2,\dots n\}$ as a nonoriented graph in which any two vertices are joined 
by an unordered edge.  Since the distributions of the random variables 
$x_{ij}$ are symmetric, we conclude that if a path $ \ \p \ $ gives
a nonzero 
contribution to (3.4) then
the following condition {\bf C2} also must hold :
$$ \ {\bf C2.} \ \ {\rm The \ \ number \ \ of \ \ occurrences \ \ of \ \ each 
\ \ edge \ \ is \ \ even.} $$
  Indeed, due to the independence of $\{ x_{ij}\} $, 
the mathematical expectation of the product factorizes as a product of 
mathematical expectations of random variables corresponding to different 
edges of the path.  Therefore if some edge appears in $\p$ odd number of 
times at least one factor in the product will be zero. 
Condition {\bf C2} is a necessary but not sufficient condition on
$ \ \p \ $ to give a non-zero contribution in (3.4). To obtain a necessary and
sufficient condition let us note that an
edge $ \ \ i_k=j, \ \ i_{k+1}= g,  \ \ k=0, \ldots, 2 \* m_p -1, 
\ \ $  contributes 
$  \ \ x_{jg} \ \ $ for  odd $k$ and    $ \ \ x_{gj} \ \ $
for even $k.\ $ Clearly the number of apperances in 
each non-zero term of (3.4)
must be even both for $ \ x_{jg} \ $ and $ \ x_{gj} \ .$ 
This leads to 
$$ \ {\bf C3.} \ \ {\rm For\   any\  edge}
 \ \{j,g\}, \ \  j,g \in \{1,2,\ldots, n\}, \ \ 
{\rm the \  number \  of \  times  \   we \  pass}$$
$$ \{j,g\} \ {\rm in \  the \  direction}  \ j \to g \ 
{\rm at \ odd \ moments \ of \ time}
 \ 2k+1, \ \ k=0,1,\ldots, m_p,$$
$$ {\rm plus \  the \  number \  of\ 
times \  we  \ pass}
 \ \{j,g\} \ \  {\rm  \ in  \ the  \ direction}  \ g \to j \  {\rm 
at  \ even}$$
$$ {\rm  moments  \ of  \ time}
 \ 2k, \ \ k=0,1,\ldots \  {\rm  must  \ be  \ even. }$$

Let us now consider the r.h.s. of (3.3). We can write
\begin{equation}
{\bf E}
 \ \Trace\ M^{2\*m_p}_n=  \sum_{\p} {\bf E}\ y_{i_0,i_1}
\* y_{i_1, i_2}\* y_{i_2, i_3}\* y_{i_3,i_4}
\dots y_{i_{2\* m_p -2}, i_{2 \* m_p -1}}
\* y_{i_{2 \* m_p-1, i_0}},
\end{equation}
where the sum  again is  over all closed paths 
$\p =\{i_0,i_1,\ldots ,i_{2\* m_p -1,i_0}\}$, with a distinguished origin, 
in the set $\{1,2,\dots n\}$. Since $ \ M_n \ $ is a square $ n \times n \ $ 
real symmetric matrix conditions {\bf C1} and {\bf C3} are no longer needed.
In particular the necessary and sufficient condition on a path $ \ \p \ $
to give a non-zero contribution to (3.5) is {\bf C2}. It does not matter in 
which direction we pass an edge $ \ \{ jg\} \ $, because both steps 
$ \ j \to g \ $ 
and $ \ g \to j \ $ give us $ \ y_{jg}=y_{gj} \ $. Using the inequalities
$ \ \  {\bf E} \ x_{jg}^{2\*r}\*  {\bf E } x_{gj}^{2\* q}
\leq {\bf E} \* y_{jg}^{2\*r +2\*q} \ $ we show 
that each term in (3.4) is not greater than the 
corresponding term in (3.5) and, therefore, obtain (3.3). (3.1) (in the case
$ \ k=1 \ $)
then immediately follows from Theorem 3 of [30] (the matrix $\ A_n \ $ 
considered there differs from $ \ M_n \ $ by a factor $ \ \frac{1}{2\*\sqrt{n}}
\ ).$  In the  general case the proof of (3.1)-(3.2) is essentially identical 
to the one given in [30]. In particular, part b) of Lemma 3
follows from the fact that the l.h.s. at (3.2) is given by a subsum over
paths that, in addition to {\bf C1-C3} have at least one edge appeared four 
times or more. As we showed in [30] the contribution of such paths tends to 
zero as $ \ n \to \infty \ \ .$ Lemma 3 is proven.
\medskip

\noindent{\bf Remark 7}

If the condition $ \ \ n-p=O(p^{1/3}) \ $ in Lemma 3 and Theorem 2 is not 
satisfied the machinery
from [28]-[30] does not work, essentially for the following reason:
when we decide which vertice to choose during the moment of self-intersection
(as explained in  section 4 of [29]) the number of choices for odd moments 
of time is 
smaller because of the constrain {\bf C1}. If we now use the same bound as 
for the even 
moments of time (the one similar to the bound at the bottom of p.725 of [29])
the estimate becomes rough when $ \  n-p \ $ is much greater then 
$ \ p^{1/3} \ .$ Therefore new combinatorial ideas are needed.
\medskip

As  corollaries of Lemma 3 we obtain 
\medskip

\noindent{\bf Corollary 2}
{\it There exist  $ \ const >0 \ $ 
such that for any $ \ s=o(p^{1/3}) \ $
$$ {\bf P } ( \lambda_1 (A_p) > \mu_{n,p} + \sigma_{n,p}\* s ) < const  \*
\exp (- const \* s) $$}
\medskip

\noindent{\bf Corollary 3}
{\it
\begin{align}
&\int_{(-\infty, \mu_{n,p} +\sigma_{n,p} \* p^{1/6}]^k}
\* \exp( \sum_{j=1,\ldots, k} t_j\*s_j) \*
\bar \rho_k^{(p)}(s_1,\ldots, s_k) ds_1 \ldots ds_k  \ \ \to \\
&\int_{R^k} \* \exp( \sum_{j=1,\ldots, k} t_j\*s_j) \*
\rho_k(s_1,\ldots, s_k) ds_1 \ldots ds_k  
\end{align}
for any $ t_1 >0, \ \ldots, t_k >0 \ $ as $ \ n \to \infty$,
where $$ \ \bar \rho_k^{(p)} (s_1, \ldots, s_k) = (\sigma_{n,p})^k \*
\rho_k^{(p)} (\mu_{n,p} +\sigma_{n,p} \*s_1, \ldots, 
\mu_{n,p} +\sigma_{n,p} \*s_k) $$ is the rescaled $k$-point correlation 
function and
$ \ \rho_k(s_1,\ldots,s_k) \ $ is defined in Section 2, Remark 6 by the r.h.s.
of (2.18)-(2.21). }

\medskip

To prove Corollary 2 we use the Chebyshev  inequality
$$ {\bf P} ( \lambda_1 (A_p) > \mu_{n,p} + \sigma_{n,p}\* s ) \leq \frac{
{\bf E} \lambda_1(A_p)^{\sigma_{n,p}}}{ (\mu_{n,p} +\sigma_{n,p} \* s)^
{p^{2/3}}} 
\leq  \frac{ {\bf E} \Trace \* A_p^{\sigma_{n,p}}}
{ (\mu_{n,p} +\sigma_{n,p} \* s)^{p^{2/3}}} $$
and Lemma 3.
As a result of Corollary 2 we obtain that with probability $ \ O(\exp(-const \*
p^{1/6})) \ $  the largest eigenvalue is not greater than
$ \ \mu_{n,p} + \sigma_{n,p}\*p^{1/6} \ .$ Therefore, it is enough to study
only the eigenvalues in $ \ (-\infty, \ \mu_{n,p} +\sigma_{n,p} \* p^{1/6}]
 \ $ (with very high probability there are no eigenvalues outside). To prove 
Corollary 3 we first note that Lemma 3 implies

\begin{align}
&\int_{(-\infty, \mu_{n,p} +\sigma_{n,p} \* p^{1/6}]^k}
\* \exp( \sum_{j=1,\ldots, k} t_j\*s_j) \*
\bar \rho_k^{(p)}(s_1,\ldots, s_k) ds_1 \ldots ds_k  \leq \\
&\frac{const^k }{\prod^k_{i=1}
t_i^{3k/2}}\*\exp (const \cdot\sum^k_{i=1}t^3_i),
\end{align}
with some  $ \ \ const  >0 .\  \ $
To see this we write
$$ {\bf E} \ \sum^* \ \prod_{j=1}^{j=k} \ 
\exp \bigl ( t_j \*(\lambda_{i_j}-\mu_{n,p})/\sigma_{n,p} \bigr ) \leq 
{\bf E} \ \sum^* \ 
\prod_{j=1}^{j=k} \ (\lambda_{i_j}/\mu_{n,p})^{[t_j \* \mu_{n,p}/\sigma_{n,p}]}
\bigl ( 1+ o(1) \bigr ) $$
$$ \leq {\bf E} \ \prod_1^k \ \Trace \ A_p^{[t_j \*\mu_{n,p}/\sigma_{n,p}]} \
 \mu_{n,p}^{- (\sum_1^k \ t_j)\*\mu_{n,p}/\sigma_{n,p}} \* 
\bigl ( 1+o(1) \bigr ) $$
where the sum $ \ \sum^* \ $ is over all k-tuples of non-coinciding indices
$ \ (i_1, i_2, \ldots, i_k),$ \\ 
$ \ \ 1 \leq i_j \leq p, \ j=1,\ldots, k, \ \ $
such that $ \ \lambda_{i_j} < \mu_{n,p} +\sigma_{n,p} \* p^{1/6},  \ \ 
j=1,\ldots,k, \  \ \ $ and apply Lemma 3, a). Part b) of Lemma 3 implies that
the differences between left hand sides of (3.8) for different ensembles of 
random matrices (i)-(iv)
( (i')-(iv')) tend to 0. Finally we note that in the Gaussian case the l.h.s.
of (3.8) converges, which in turn implies the convergence for arbitrary 
ensemble of sample covariance matrices. For the details we refer the reader to
the analogous arguments in [30]. Corollary 3 is proven.
\medskip
Theorem 2 now follows from Lemma 2, part d) and Corollary 3.
\medskip

\section{ Proof of Theorem 3}
In order to estimate the r.h.s. of (3.4)
we assume  some familiarity of the reader with the combinatorial 
machinery developed in [28]-[30].
In particular we refer the reader to [28] ( Section 2, Definitions 1-2) or
[29] ( Section 4, Definitions 1-4) how we defined  a)  
marked and unmarked instants, b) a partition of all verices into the 
classes $ \ \ \n_0, \n_1, \ldots, \n_m \ \ $ and c) paths of  the type 
$ \ (n_0,n_1, \ldots, n_m) \  \ $, where $ \ \ 
 \sum_0^m \ n_k =n, \ \ \sum_0^m \ k\*n_k =m 
\ \ \ $(for simplicity we  omit a subindex $p$ in $m_p$ throughout this 
section). Let us 
first estimate  a subsum of (3.4) over the paths of some fixed 
type $ \ (n_0, n_1, \ldots, n_m) \ .$
Essentially repeating the arguments from [28]-[29] we can bound it from above 
by
\begin{equation}
 \* p^{n_1+1}\* \frac{(n-n_1)!}{n_0!\*n_1!\dots n_m!}
\*\frac{m!}{\prod^{m}_{k=2}(k!)^{n_k}}\*
\prod^{m}_{k=2}(\const \* k)^{k\cdot n_k} \ 
 \sum_{X\in \Omega_m} 
\ (n/p)^{\#(X)},
\end{equation}
where  the sum $ \ \ \sum_{X\in \Omega_m} \ \ $
is over all possible trajectories $$ X=\{ x(t)\geq 0,  \ 
 \ x(t+1)-x(t)=-1, +1, \ \ t=0,\ldots, 2m-1, x(o)=x(2m)=0 \}$$
and $ \ \#(X)= \# (t : x(t+1)-x(t)=+1, \ t=2k, \ k=0,\ldots, m-1) \ \ .$

The only differences between the estimates 
(4.1) in this paper and (4.4) and (4.27) in [29]
are \\
a) the number of ways we can choose the vertices from $ \n_1 \ $
is estimated from above by $ \ \ p^{n_1} \ (n/p)^{\#(X)} /n_1! \ \ 
$ not by $ \ n\*(n-1)\cdots (n-n_1+1) \ /n_1! \ ,$
because of the restriction
{\bf C1 } from the last section ,\\
b) we have in (4.1) the factor $ \ \ (const\* 2)^{2\*n_2} \ \ $ instead
of $ \ \ 3^r \ \ $ in (4.27) of [29], which is perfectly fine since
$ \ \ r\leq n_2 \ \ $ (by $ \ r \ $ we denoted in [29] the number of so-called
``non-closed'' verices from $ \ \n_2 \ ),$ and \\
c) there is no factor $ \ \ \frac{1}{n^m} \ \ $ in (4.1) because of the 
different normalization.
Let us denote by
$ \ \ g_m(y)=\sum_{X\in\Omega_m} \ y^{\#(X)} \ \ $(observe that $ \  g_m(1)=
|\Omega_m|=\frac{2m!}{m!\*(m+1)!} \ \ $ are just Catalan numbers).\\
Consider the generating function
$ \ \ G(z,y)= \sum_{m=0}^{\infty} \ g_m(y)\*z^m, \ \ g_0(y)=1 \ \ .$
It is not difficult to see  (by representing $ \ g_m(y) \ $ as a sum over the 
first 
instants of the  return of the trajectory to the origin) that
\begin{align}
\begin{split}
& G(z,y)=1+y\*z\* G'(z,y)\*G(z,y)\\
& G'(z,y)=1+z\* G'(z,y)\*G(z,y),
\end{split}
\end{align}
where
\begin{align}
\begin{split}
&G'(z,y)=\sum_{m=0}^{\infty} \ {g'}_m(y)\* z^m, \ \ \ \ 
 {g'}_m(y) = \sum_{X\in\Omega_m} \ y^{\#'(X)} \ \ \ \ {\rm and}\\
&\#'(X)= \# (t : x(t+1)-x(t)=+1, \ t=2k+1, \ k=0,\ldots, m-1).
\end{split}
\end{align}
Solving (4.2) we obtain 
\begin{equation}
G(z,y)=\frac{-y\*z +z +1 -\sqrt{((y-1)\*z -1)^2 -4z}}{2z}=
\frac{-y\*z +z +1 -(y-1)\sqrt{(z-z_1)(z-z_2)}}{2z},
\end{equation}
where $ \ \ z_1=1/(\sqrt{y}+1)^2, \ \ z_2=1/(\sqrt{y}-1)^2, \ \ $ and we take 
the branch of $ \ \ \sqrt{(z-z_1)\*(z-z_2)}, \ \ $ analytic everywhere outside
$ \ [z_1, z_2], \ $  such that
that $ \ \ \sqrt{(0-z_1)\*(0-z_2)}=1/(y-1) \ \ .$ Therefore
\begin{equation}
g_m(y)= -\frac{y-1}{4\pi \*i}\  
\oint_{|z|=z_1-\epsilon} \  \frac{\sqrt{(z-z_1)(z-z_2)}}{z^{m+2}},\ \ 
\ m\geq 1,
\end{equation}
where the integration is counter-clockwise. An exercise in complex analysis
gives us
\begin{align}
\begin{split}
\oint_{|z|=z_1-\epsilon} \  \frac{\sqrt{(z-z_1)(z-z_2)}}{z^{m+2}} &=
-2i \*\frac{\sqrt{z_2-z_1}}{z_1^{m+1/2}\* m^{3/2}} \ \int_0^{\infty} \ \sqrt{t}
\* \exp(-t) \* dt \ (1+o(1))\\
& = \frac{ 2i\* \sqrt{\pi} \* y^{1/4} \* (\sqrt{y}+1)}{(y-1)} \* 
\frac{(\sqrt{y}+1)^m}{m^{3/2}} \ (1+o(1)).
\end{split}
\end{align}
Therefore 
\begin{equation}
 g_m(y)=\frac{ y^{1/4} \* (\sqrt{y}+1)}{2\*\sqrt{\pi}} \* 
\frac{(\sqrt{y}+1)^m}{m^{3/2}} \ (1+o(1)),
\end{equation}
and the subsum of (3.4) over the paths of the type $ \ \ (n_0, n_1, 
\ldots, n_m) \ \ $ is bounded from above by
\begin{align}
\begin{split}
&\frac{ (n/p)^{1/4} \* (\sqrt{n/p}+1)}{2\*\sqrt{\pi}} 
\  p^{n_1+1}\* \frac{(n-n_1)!}{n_0!\*n_1!\dots n_m!}
\*\frac{m!}{\prod^{m}_{k=2}(k!)^{n_k}}\\
&\prod^{m}_{k=2}(\const \* k)^{k\cdot n_k} 
\ \frac{(\sqrt{n/p}+1)^m}{m^{3/2}} \ (1+o(1)) \leq\\
& \frac{ (n/p)^{1/4} \* (\sqrt{n/p}+1)}{2\*\sqrt{\pi}} 
\ \frac{p \ \mu_{n,p}^m}{m^{3/2}} \ \frac{1}{p^{m-n_1}} \ 
\frac{(n-n_1)!}{n_0!\*n_1!\dots n_m!}
\*\frac{m!}{\prod^{m}_{k=2}(k!)^{n_k}}
\ \prod^{m}_{k=2}(\const \* k)^{k\cdot n_k} \ (1+o(1))
\end{split}
\end{align}
(the constant $ \ const \ $ may have changed).
Using the inequality
$ \ \ m! <n_1!\* m^{m-n_1} \ \ $ and $ \ \ \sum_{k=1}^m k\*n_k=m, \ \ 
\sum_{k=1}^{m}=n-n_0 \ \ $
we obtain 
\begin{align}
\begin{split}
(4.8) &\leq 
\frac{ (n/p)^{1/4} \* (\sqrt{n/p}+1)}{2\*\sqrt{\pi}} 
  \ \frac{p \ \mu_{n,p}^m}{m^{3/2}} \
n^{ -\sum_2^k k\*n_k} \ n^{\sum_2^m n_k} \ m^{\sum_2^m k\*n_k} \ 
\prod_{k=2}^m \ \frac{(const^k)^{n_k}}{n_k!}\\
&\leq \frac{ (n/p)^{1/4} \* (\sqrt{n/p}+1)}{2\*\sqrt{\pi}} 
  \ \frac{p \ \mu_{n,p}^m}{m^{3/2}} \
\bigl ( \prod_2^m \ \frac{1}{n_k!} \ \bigl( \frac{ (const \ m)^k}{n^{k-1}}
\bigr )^{n_k} \bigr )
\end{split}
\end{align}
Now we can  estimate the sum  of (4.9) over  $ \ (n_o, n_1, \ldots, n_m), $\\ 
$  0<\sum_{k=2}^{m}
\ k\*n_k \leq m \ \ $ as
\begin{equation}
\frac{ (n/p)^{1/4} \* (\sqrt{n/p}+1)}{2\*\sqrt{\pi}} 
  \ \frac{p \ \mu_{n,p}^m}{m^{3/2}} 
\  \bigl ( \exp( \sum_{k=2}^{m} \frac{ (const \ m)^k}{n^{k-1}})-1\bigr )
\end{equation}

Since for $ \ \ m=o(p^{1/2}) $
\begin{equation}
\sum_{k=2}^{m} \frac{ (const \  m )^k}{n^{k-1}} =O(m^2/n)=o(1) 
\end{equation}
we see that the subsum of (3.4) over $ \ \p \ $ with
$ \ \ \sum_{k=2}^m n_k >0 \ \ $ is $o( \frac{p \ \mu_{n,p}^m}{m^{3/2}}).$
Finally we note that the subsum over the paths of the type $ \ \ (n-m,m,0,0,
\ldots, 0) \ \ $
is 
\begin{equation}
\frac{ (n/p)^{1/4} \* (\sqrt{n/p}+1)}{2\*\sqrt{\pi}} 
  \ \frac{p \ \mu_{n,p}^m}{m^{3/2}} \ (1+o(1)),
\end{equation}
because  for such paths we can choose the vertices from $ \ \n_1 \ $ exactly
in $ \ p^m \* (n/p)^{\#(X)} \*(1+o(1)) \ $ different ways ( if
$ \ m=o(p^{1/2}) \ $), and the first point of 
a path in 
$ \ p \ $ different ways. Combining (4.11) and (4.12) we prove the first
part  of Theorem 3. To prove part b)  we observe that if 
$ \ \ m =O(p^{1/2}), \ \  $
 the l.h.s. of (4.11)
is still $ \ \ O(m^2/n), \ \ $ which together with (4.10) and (4.12) finishes 
the proof.
Theorem 3 is proven.\\
To derive  Corollary 1  from Theorem 3 we apply   Chebyshev's inequality
(similarly to the proof of Corollary 2 in Section 3) and  Borel-Cantelli 
lemma.

\def\am{{\it Ann. of Math.} }
\def\ap{{\it Ann. Probab.} }
\def\temf{{\it Teor. Mat. Fiz.} }
\def\jmp{{\it J. Math. Phys.} }
\def\cmp{{\it Commun. Math. Phys.} }
\def\jsp{{\it J. Stat. Phys.} }
\def\cpam{{\it Comm. . Pure  Appl. Math.} }


\begin{thebibliography}{99}


\bibitem[1] {} D. Aldous, P. Diaconis, Longest increasing subsequences: 
From patience sorting to the Baik-Deift-Johansson theorem, {\it 
Bull. Amer. Math. Society (N.S.)}  {\bf 36} 413--432 (1999).


\bibitem[2] {} Z.D. Bai, Y.Q. Yin,
Convergence of the semicircle law,
\ap {\bf 16},  863--875 (1988).

\bibitem[3]{} Z.D. Bai, J. W. Silverstein, Y.Q. Yin, A note on the largest 
eigenvalue of a large dimensional
sample covariance matrix, {\it J. Mult. Anal.} {\bf 26},  166--168 (1988). 

\bibitem[4]{} Z.D. Bai, Convergence rate of expected spectral distributions
of large random matrices.  II.  Sample covariance matrices, 
\ap {\bf 21},  649--672
(1993).

\bibitem[5] {} Z.D. Bai, J.W. Silverstein, No eigenvalues outside the 
support of the limiting spectral distribution of large-dimensional sample 
covariance matrices, {\ap} {\bf 26}, No.1, 316--345, 1998.


\bibitem[6] {} Z.D. Bai, Methodologies in spectral analysis of large 
dimensional random matrices. A review, {\it Statistica Sinica} {\bf 9}, No. 3,
611--677 (1999).


\bibitem[7]{} P. Bleher, A. Its, Semiclassical asymptotics of orthogonal
polynomials, Riemann-Hilbert problem, and universality in the matrix model,
\am, {\bf 150}, No.1, 185-266, (2000).

\bibitem[8]{} E. Br\'ezin, A. Zee, Universality of the correlations between
eigenvalues of large random matrices, {\it Nuclear Phys. B.} {\bf 402},
 613--627 (1993).


\bibitem[9]{} B.V. Bronk,
Exponential ensemble for random matrices, {\it J. Math. Phys.} {\bf 6},
No.2, 228--237  (1965).


\bibitem[10] {} P. Deift, Integrable Systems and Combinatorial Theory,
{\it Notices of the AMS} {\bf 47}, No. 6, 631--640
  (2000).


\bibitem[11]{} P. Deift, T. Kriecherbauer, K.T.-R. McLaughlin, 
S. Venakides,
X. Zhou, Uniform asymptotics for polynomials orthogonal with respect to
varying exponential weights and applications to universality questions in
random matrix theory,  \cpam , {\bf 52}, No.11, 1335-1425, (1999).


\bibitem[12]{} A. Edelman, The distribution and moments of the smallest
eigenvalue of a random matrix of Wishart type,
{\it Linear  Algebra and  its Applications} {\bf 159}, 55--80 (1991).


\bibitem[13]{}P.J. Forrester, The spectral edge of random matrix ensembles,
{\it Nucl. Phys. B.} {\bf 402},  709--728 (1994).


\bibitem[14]{} P.J. Forrester, T.Nagao, G.Honner, Correlations for the 
orthogonal-unitary and symplectic-unitary transitions at the hard and soft 
edges, 
{\it Nucl. Phys. B.} {\bf 553}, 601-643 (1999).


\bibitem[15]{} S. Geman, A limit theorem for the norm of random matrices,
{\it Ann. Probab.} {\bf
8}, 252--261 (1980).


\bibitem[16] {} A. Guionnet, Large deviation upper bounds and central limit
theorems for band matrices and non-commutative functionals of Gaussian large
random matrices, preprint (2000).

\bibitem[17]{} U. Grenader, J. Silverstein, Spectral analysis of networks 
with random topologies, {\it SIAM J. Appl. Math. } {\bf 32},  
499--519 (1977). 


\bibitem[18]{} A.T. James,
Distribution of matrices variates and latent roots derived from normal samples
, {\it Ann. Math. Statist.} {\bf 35},
475-501--237  (1965).


\bibitem[19]{}K. Johansson, Shape fluctuations and random matrices,
\cmp, {\bf 209}, 437-476,  (2000).


\bibitem[20]{} K. Johansson, Universality of local spacing distributions
in certain Hermitian Wigner matrices,  \cmp  {\bf 215}, No. 3, 683-705 (2001).

\bibitem[21] {} I.M. Johnstone, On the distribution of the largest principal
component, {\it Ann. Stat.}, {\bf 29}, No. 2, (2001).

\bibitem[22] {} E. Kanzieper, V. Frelikher, Spectra of large matrices:
a method of study, {\rm in Diffusive Waves in Complex Media}, Kluwer, Dodrecht,
165--211, 1999.

\bibitem[23] {} V. Marchenko, L.Pastur, The eigenvalue distribution in some 
ensembles of random matrices, {\it Math USSR Sbornik}, {\bf 1}, 
457--483 (1967).

\bibitem[24] {} M.L.Mehta, Random Matrices, 2nd edition, Academic Press, 
New York, 1991.


\bibitem[25] {} R.J. Muirhead, Aspects of Multivariate Statistical Theory,
Wiley, New York, 1982.


\bibitem[26]{} L.A. Pastur, M. Shcherbina, Universality of the local
eigenvalue statistics for a class of unitary invariant random matrix 
ensembles, \jsp {\bf 86},  109--147 (1997).


\bibitem[27]{} J. Silverstein, On the weak limit of the largest 
 eigenvalue of a large dimensional
sample covariance matrix, {\it J. Mult. Anal.} {\bf 30},  307--311 (1989). 

\bibitem[28]{} Ya. Sinai, A. Soshnikov, Central limit theorem for traces
of large random symmetric matrices, {\it Bol. Soc. Brasil. Mat.},
{\bf 29}, No. 1,  1--24 (1998).

\bibitem[29]{} Ya. Sinai, A. Soshnikov, A refinement of Wigner's semicircle
law in a neighborhood of the spectrum edge for random symmetric matrices,
{\it Functional Anal. Appl.} {\bf 32}, No. 2, (1998).

\bibitem[30] {} A. Soshnikov, Universality at the edge of the spectrum in 
Wigner random matrices,
\cmp, {\bf 207}, 697--733 (1999).

\bibitem[31] {} A. Soshnikov, Determinantal random point fields,
{\it Russian Math Surveys}, {\bf 55}, No.5, 923-975, (2000).

\bibitem[32] {} G. Szego, Orthogonal Polynomials, 3rd edition, 
American Mathematical Society, 1967.


\bibitem[33]{} C.A. Tracy, H. Widom, Level-spacing distribution and Airy
kernel, \cmp {\bf 159},  151--174 (1994).

\bibitem[34]{} C.A. Tracy, H. Widom, On orthogonal and symplectic matrix
ensembles, \cmp {\bf 177},  727--754 (1996).

\bibitem[35]{} C.A. Tracy, H. Widom, Correlation functions, cluster functions,
and spacing distribution for random matrices, \jsp {\bf
92}, No. 5/6, (1998).


\bibitem[36]{} Y.Q. Yin, Z.D. Bai, P.R. Krishnaiah, On the limit of the 
largest eigenvalue of a large dimensional
sample covariance matrix, {\it Probab, Theory  Related Fields } {\bf 78},  
509--521 (1988). 


\bibitem[37]{} K.W. Wachter, The strong limits of random matrix spectra
for sample matrices of independent elements, {\it Ann. Probab.} {\bf
6}, No. 1,  1--18 (1978).


\bibitem[38]{} H. Widom, On the relation between orthogonal, symplectic
and unitary matrix ensembles, \jsp  {\bf 94}, No 3/4, (1999).


\bibitem[39] {} S.Wilks, Mathematical Statistics, Wiley, New York, 1967.






\end{thebibliography}
\end{document}